\documentclass[a4paper,10pt]{article}
\usepackage{amssymb}
\usepackage{amsmath}
\usepackage{nicefrac}
\usepackage{graphicx}
\usepackage{dcolumn}
\usepackage{bm}
\usepackage{comment}
\usepackage{multirow}
\usepackage{cite}
\usepackage{xcolor}
\usepackage{url}
\usepackage{mathrsfs}
\usepackage[normalem]{ulem}
\usepackage{bold-extra}
\usepackage{hyperref}
\usepackage{amsthm} 
\usepackage[utf8]{inputenc}
\usepackage[T1]{fontenc}
\usepackage[english]{babel}
\usepackage{algorithm}
\usepackage{algpseudocode}
\usepackage{booktabs}

\newcommand{\N}{\mathbb{N}}

\newcommand{\Sg}{\mathcal{S}}
\newcommand{\Gap}{\mathcal{F}}

\theoremstyle{definition}

\begin{document}

\huge

\begin{center}
Steganographic information hiding via symmetric numerical semigroups
\end{center}

\vspace{0.5cm}

\large

\begin{center}
Jean-Christophe Pain$^{a,b,}$\footnote{jean-christophe.pain@cea.fr}
\end{center}

\normalsize

\begin{center}
\it $^a$CEA, DAM, DIF, F-91297 Arpajon, France\\
\it $^b$Universit\'e Paris-Saclay, CEA, Laboratoire Mati\`ere en Conditions Extr\^emes,\\
\it 91680 Bruy\`eres-le-Ch\^atel, France\\
\end{center}

\begin{abstract}
We introduce a steganographic information hiding scheme based on structural properties of numerical semigroups arising from the Frobenius coin problem. Instead of encoding data through representable integers, the proposed protocol embeds information into the gap structure of carefully chosen symmetric numerical semigroups. Symmetry guarantees a balanced gap density, ensuring that encoded values are statistically indistinguishable from uniform numerical noise to an observer lacking the private generating set. The security of the scheme relies on the assumed average-case hardness of numerical semigroup membership inference for hidden generators, offering a novel number-theoretic primitive for covert communication and post-quantum resilient information hiding.
\end{abstract}

\section{Introduction}\label{sec1}

The Frobenius ``coin problem'' (also known as the ``money changing problem'') originates from a deceptively simple question in additive number theory: given a set of coprime positive integers $A = \{a_1, a_2, \dots, a_n\}$, what is the largest integer $F(A)$ that cannot be expressed as a non-negative integer linear combination of its elements? While the determination of this value, known as the Frobenius number, was a primary focus for 19th-century mathematicians like Sylvester \cite{Sylvester1884} and Frobenius \cite{Frobenius1894}, its modern study is fundamentally intertwined with the theory of numerical semigroups \cite{Kaplan2017}.

In this algebraic context, a numerical semigroup $S$ is a commutative monoid under addition, comprising all non-negative integers except for a finite set of omitted values, or ``gaps'' $\mathbb{N} \setminus S$. The Frobenius number $F(S)$ thus represents the maximum of this set of gaps, marking the threshold, often referred to as the conductor minus one, beyond which all integers are representable. During the second half of the twentieth century, interest in these structures resurfaced significantly, driven by their deep applications in algebraic geometry and commutative algebra \cite{Rosales2009}.

The computational complexity of the problem reflects its structural depth; finding the Frobenius number without restrictions on the generating set is known to be NP-hard \cite{Ramirez1996,Ramirez2005}. No polynomial-time algorithm is currently known to recover the generating set from observed gaps without additional structural information, highlighting the potential hardness for information hiding applications. Consequently, closed-form solutions are highly prized and typically restricted to specific cases. Beyond Sylvester's classical result for two generators ($n=2$, see for instance an inductive proof in Ref. \cite{Kapetanakis2021} or a very short proof in Appendix \ref{appA}), efficient algorithms for the three-generator case were developed by Greenberg \cite{Greenberg1988,Greenberg1988b}, while Roberts \cite{Roberts1956,Roberts1956b} and Lewin \cite{Lewin1975,Lewin1975b} provided foundational formulas for arithmetic and almost arithmetic sequences, respectively. It is worth mentioning that a generalization of Sylvester's formula was recently proposed by Binner \cite{Binner2021}. General properties of Frobenius number are given in Appendix \ref{appB}.

For a set in an arithmetic progression defined by $a, d, w$ with $\gcd(a, d) = 1$, the Frobenius number is given by:
\begin{equation}
F(a,a+d,a+2d,\dots ,a+wd)=\left(\left\lfloor \frac {a-2}{w}\right\rfloor \right)a+d(a-1).
\end{equation}
Interestingly, the stability of this number is remarkable. For instance, when $a>w^{2}-3w+1$, the removal of specific intermediate elements from the progression does not alter the Frobenius number. This phenomenon shifts the analytical focus from merely calculating $F(A)$ to exploiting the distribution of the integers it misses: the gaps.

The paper is organized as follows. Section \ref{sec2} establishes the mathematical framework, with particular emphasis on symmetric numerical semigroups and their gap structure. In Section \ref{sec3}, we present the steganographic protocol, including key generation and the encoding strategy based on modular gap partitioning. In the same section, the decryption algorithm and its efficient implementation via residue graph techniques are also presented. Section \ref{sec4} discusses the security properties of the scheme, clarifies the underlying hardness assumptions, and outlines its limitations. Finally, the appendix \ref{appC} provides implementation details and  a prototype of the proposed construction.

\section{Mathematical preliminaries}\label{sec2}

\subsection{Numerical semigroups and gaps}\label{subsec21}

This section recalls the main notions from numerical semigroup theory \cite{Howie1985} required throughout the paper and introduces the structural tools that underlie both the encoding and decoding procedures. A numerical semigroup $\Sg$ is a submonoid of $(\N,+)$ satisfying the following properties:
\begin{itemize}
    \item $0 \in \Sg$,
    \item $\Sg + \Sg \subset \Sg$,
    \item $\N \setminus \Sg$ is finite.
\end{itemize}
Equivalently, a numerical semigroup is a finitely generated additive submonoid of $\N$ whose generators have greatest common divisor equal to one. If $A=\{a_1,\dots,a_n\}$ is such a generating set, we write $\Sg=\langle A\rangle$.

The complement $\Gap(\Sg)=\N\setminus\Sg$ is called the set of gaps. Its cardinality
\[
g(\Sg)=|\Gap(\Sg)|
\]
is the genus of the semigroup, while the largest gap is the Frobenius number, denoted by $F(\Sg)$. Every integer greater than $F(\Sg)$ belongs to $\Sg$. The distribution of gaps below the Frobenius number encodes much of the arithmetic structure of the semigroup and forms the central object exploited by the proposed protocol.

Let $m=\min(\Sg\setminus\{0\})$ denote the multiplicity of $\Sg$. The Ap\'ery set of $\Sg$ with respect to $m$ is defined as
\[
\mathrm{Ap}(\Sg,m)=\{s\in\Sg \mid s-m \notin \Sg\}.
\]
The Ap\'ery set contains exactly $m$ elements, one in each residue class modulo $m$, and provides a complete description of the semigroup: an integer $x$ belongs to $\Sg$ if and only if
\[
x \geq \mathrm{Ap}(\Sg,m)[x \bmod m].
\]
This characterization underlies the decoding algorithm presented later. In practice, the Ap\'ery set can be computed via a shortest-path algorithm on a directed graph whose vertices correspond to residue classes modulo $m$ and whose edges are weighted by the generators of $\Sg$ \cite{Herzog1970,Paillier1999}.

\subsection{Symmetric semigroups and density}\label{subsec22}

A numerical semigroup $\Sg$ is said to be symmetric if for every integer $z$,
\[
z \in \Sg \iff F(\Sg) - z \notin \Sg.
\]
Symmetry induces a perfect duality between representable integers and gaps below the Frobenius number. In particular, symmetric semigroups satisfy the identity
\[
g(\Sg) = \frac{F(\Sg)+1}{2},
\]
which implies that exactly half of the integers in the interval $[0,F(\Sg)]$ are gaps. This symmetry implies a perfect pairing between elements of $\Sg$ and its gaps inside the interval $[0,F(\Sg)]$. In other words, this ensures that gaps occur with asymptotic density $1/2$. This balanced gap density plays a key role in ensuring steganographic indistinguishability, as it prevents the introduction of detectable statistical bias.

A generating sequence $A=(a_1,\dots,a_n)$ is called telescopic if, for all $i>1$,
\[
\frac{a_i}{\gcd(a_1,\dots,a_{i-1})} \in \langle a_1,\dots,a_{i-1}\rangle.
\]
Numerical semigroups generated by telescopic sequences are known to be symmetric and to admit particularly simple Ap\'ery sets. These structural properties ensure both efficient membership testing for legitimate receivers and a well-controlled gap structure, making telescopic semigroups especially suitable for cryptographic and steganographic constructions.

Finally, for symmetric numerical semigroups, the asymptotic density of gaps is $1/2$. When restricting attention to residue classes modulo a fixed integer $M$, the gap distribution remains approximately uniform across classes, provided $M$ is small relative to the generators. This observation justifies the modular gap partitioning strategy employed in the encoding phase: each gap can be viewed as carrying $\log_2 M$ bits of information, yielding a flexible trade-off between embedding rate and numerical range.

\section{The protocol}\label{sec3}

The proposed scheme encodes information by exploiting the arithmetic structure of gap positions within a secret numerical semigroup. The legitimate receiver possesses the generating set defining the semigroup, while an external observer only sees a stream of integers without structural context.

\subsection{Key generation and telescopic generating sets}\label{subsec31}

The private key consists of a generating set $A = (a_1,\dots,a_n)$ defining the numerical semigroup $\Sg = \langle A \rangle$. To ensure symmetry and algorithmic tractability for the legitimate receiver, we restrict to telescopic generating sets. As mentioned in the preceding section, telescopic semigroups are known to be symmetric, which directly enforces the balanced gap density required for statistical indistinguishability. Moreover, this structure enables efficient membership testing via residue graph techniques.

\subsection{Encoding via modular gap partitioning}\label{subsec32}

Rather than mapping individual bits directly to gap membership, we partition the gap set using residue classes. For a fixed modulus $M=16$, a 4-bit nibble $V\in\{0,\dots,15\}$ is encoded by selecting a value
\[
x \in \Gap(\Sg) \quad \text{such that} \quad x \equiv V \pmod{M}.
\]

This approach spreads encoded values across the numerical range and avoids low-entropy encodings. Provided the semigroup is symmetric, each residue class contains asymptotically the same number of gaps, preserving uniformity.

\subsection{Decryption algorithm}\label{subsec33}

Decoding relies on the ability of the legitimate receiver to efficiently test membership in $\Sg$. This is achieved through the construction of a residue graph modulo the smallest generator $a_1$. For each residue class modulo $a_1$, the algorithm computes the minimal representable value using shortest-path techniques. Any integer smaller than this threshold is necessarily a gap.

\section{Security analysis}\label{sec4}

This section analyzes the security properties of the proposed construction from a steganographic perspective. The scheme is not intended to provide semantic encryption on its own, but rather to serve as an information-hiding layer whose security relies on the difficulty of inferring structural properties of numerical semigroups from unlabeled numerical data.

\subsection{Adversarial model}\label{subsec41}

We consider a passive adversary observing a stream of integers $(x_i)_{i\geq 1}$ generated by the encoder. The adversary does not possess the private generating set $A = (a_1,\dots,a_n)$ defining the numerical semigroup $\Sg = \langle A \rangle$, nor any auxiliary information about its structure.

The adversary's goal is to distinguish values encoding information from uniformly distributed numerical noise, or to recover the embedded payload with non-negligible advantage. The adversary is assumed to be computationally bounded and may perform arbitrary polynomial-time precomputation, but does not have access to adaptive queries or oracle feedback.

\subsection{Hardness assumption}\label{subsec42}

The security of the scheme relies on the assumed difficulty of the so-called ``numerical semigroup membership inference problem'': given an integer $x$ and samples drawn from an unknown numerical semigroup $\Sg$, determine whether $x$ belongs to $\Sg$ or to its gap set $\Gap(\Sg)$ without knowledge of a generating basis.

While the numerical semigroup membership problem is known to be NP-hard in the worst case, we do not claim a formal worst-case-to-average-case reduction. Instead, security relies on average-case hardness for generating sets chosen to be non-superincreasing, non-redundant, and of comparable magnitude. No polynomial-time algorithm is currently known to infer the hidden generating set from unlabeled gap observations in the average case.

\subsection{Statistical indistinguishability via symmetry}\label{subsec43}

A core feature of the protocol is the use of symmetric numerical semigroups. For such semigroups, the symmetry relation
\[
z \in \Sg \iff F(\Sg) - z \notin \Sg
\]
induces a perfect pairing between representable integers and gaps within the interval $[0,F(\Sg)]$. Consequently, exactly half of the integers in this interval are gaps.

This balanced density implies that gap membership behaves approximately as an unbiased Bernoulli random variable when restricted to bounded numerical windows. As a result, encoded values are statistically indistinguishable from uniformly distributed numerical noise under standard frequency-based statistical tests, provided that the structure of $\Sg$ remains hidden.

\subsection{Resistance to lattice-based attacks}\label{subsec44}

Classical lattice reduction attacks, such as LLL (Lenstra–Lenstra–Lov\'asz) \cite{Lenstra1982} or BKZ (Block Korkine-Zolotarev) \cite{Schnorr1987,ChenNguyen2011}, are effective against knapsack-type constructions when the underlying generating set exhibits strong structural bias, such as superincreasing sequences or low-density representations.

In the present construction, the adversary is not exposed to linear combinations of the generators, but only to isolated membership instances obscured by gap selection. Moreover, the generators are chosen to avoid superincreasing behavior and to maintain bounded ratios between their magnitudes. These constraints prevent the direct construction of a lattice basis amenable to reduction techniques and significantly limit the applicability of classical lattice attacks.

\subsection{Salting and range expansion}\label{subsec45}

To further reduce structural leakage, a salting mechanism is applied to encoded values:
\[
x' = x + k \cdot \mathrm{lcm}(a_i,a_j),
\]
where $(a_i,a_j)$ is a fixed pair of generators and $k$ is a random integer.

Since $\mathrm{lcm}(a_i,a_j)$ belongs to the numerical semigroup $\Sg$, adding integer multiples of this quantity preserves the gap or non-gap status of $x$ for the legitimate receiver. For an external observer lacking knowledge of the generators, this transformation expands the observable numerical range and disrupts correlations between residues and gap positions.

We emphasize that the present analysis does not constitute a formal cryptographic proof. In particular, no reduction to a standard cryptographic hardness assumption is claimed, and adaptive or active adversarial models are not considered. The proposed scheme should therefore be viewed as a steganographic primitive or covert encoding layer, suitable for composition with conventional cryptographic encryption mechanisms.

Despite these limitations, the construction introduces a novel use of additive number theory in information hiding. By leveraging the structural properties of symmetric numerical semigroups and their gap distributions, the protocol achieves balanced statistical behavior and resistance to straightforward structural inference, suggesting new directions for the application of numerical semigroup theory in secure and covert communication.

\section{Conclusion}

This work demonstrates that numerical semigroup theory provides a promising and relatively unexplored framework for steganographic information hiding. By encoding data within the gap structure of symmetric numerical semigroups, the proposed protocol achieves statistical balance, structural obfuscation, and efficient decoding for legitimate receivers.

While the present scheme does not claim formal cryptographic security in the classical sense, it introduces a novel number-theoretic perspective on covert communication and highlights new connections between additive combinatorics and information security. Future work will investigate parameter selection, empirical indistinguishability testing, and resistance against advanced lattice-based and combinatorial attacks.

\appendix

\section{The Sylvester formula}\label{appA}

The solutions of the equation $ax+by=N$ are of the form $(x_0+bt,y_0-at)$, with $t\in\mathbb{Z}$ and $(x_0,y_0)$ a particular solution of the equation. Let us assume that $N\geq ab-a-b$, and $t$ be an integer such that $\leq y_0-at\leq a-1$, then we have
$$
(x_0+bt)a=N-(y_0-at)b>ab-a-b-(a-1)b=-a,
$$
yielding $x_0+bt>-1$, i.e., $x_0+bt\geq 0$. Thus, the equation $ax+by=N$ has positive integer solutions, and
$$
F(a,b)\geq ab-a-b.
$$
If the equation $ax+by=ab-a-b$ has no positive integer solutions, then $ab=a(x+1)+b(y+1)$. Since $\gcd(a,b)=1$, we have that $a\vert(y+1)$ and $b\vert(x+1)$, which implies that $y+1\geq a$ and $x+1\geq b$. Subsequently, 
$$
ab=a(x+1)+b(y+1)\geq 2ab
$$
which is a contradiction and thus 
$$
F(a,b)\geq ab-a-b,
$$
yielding finally the Sylvester formula $F(a,b) = ab-a-b$.

\section{Special cases and asymptotic behavior of the Frobenius number}\label{appB}

\subsection{Geometric sequences and the Frobenius number}\label{appB1}

Consider positive integers \(a, b, k\) with \(\gcd(a,b)=1\), and define the geometric sequence
\[
A_k(a,b) = \{ a^k, a^{k-1}b, a^{k-2}b^2, \dots, b^k \}.
\]
Let $\Sg = \langle A_k(a,b) \rangle$ be the numerical semigroup generated by this set. We denote $F(\Sg)$ the Frobenius number, i.e., the largest integer not representable as a non-negative linear combination of elements of $A_k(a,b)$, and $g(\Sg) = |\Gap(\Sg)|$ the genus, i.e., the number of gaps in $\Sg$. For convenience, we introduce the sum notation
\begin{equation}
\sigma_r(a,b) = a^r + a^{r-1}b + a^{r-2}b^2 + \dots + b^r = \sum_{i=0}^{r} a^{r-i} b^i,
\end{equation}
i.e., the sum of all terms in the geometric progression of exponent \(r\) \cite{Ong2008,Tripathi2008}. Using this notation, the Frobenius number of \(A_k(a,b)\) is given by the Tripathi formula \cite{Tripathi2008}:
\begin{equation}
F(A_k(a,b)) = \sigma_{k+1}(a,b) - \sigma_k(a,b) - \bigl(a^{k+1} + b^{k+1}\bigr).
\end{equation}
For example, when \(k=2\), we have
\[
\sigma_2(a,b) = a^2 + ab + b^2, \quad \sigma_3(a,b) = a^3 + a^2b + ab^2 + b^3,
\]
so that
\[
F(A_2(a,b)) = \sigma_3(a,b) - \sigma_2(a,b) - (a^3 + b^3) = a^2b + ab^2 - ab = ab(a+b-1).
\]
If the semigroup is symmetric, the genus is
\[
g(\Sg) = \frac{F(\Sg)+1}{2}.
\]
For larger \(k>2\), \(\sigma_k(a,b)\) can be used to express intermediate sums, generating functions, or simplify combinatorial derivations. It is important to note that \(\sigma_k(a,b)\) itself is not the Frobenius number or the genus, but a convenient auxiliary notation to write the formula compactly. It is worth mentioning that formulae and fast algorithms \cite{Tripathi2017} are known for three numbers.

\subsection{Asymptotic bounds and average behavior}\label{appB2}

When the number of generators is small, the complexity of the "Coin Problem" remains manageable. Significant algorithmic frameworks for handling $n$ generators were established in the late 1970s, notably by Selmer \cite{Selmer1977}. However, as $n$ grows, hand calculations are often tedious. In this context, research has shifted toward sharp asymptotic bounds. Davison established a lower bound for the modified Frobenius number $f(a_{1},a_{2},a_{3}) = g(a_{1},a_{2},a_{3}) + a_{1}+a_{2}+a_{3}$, showing that:
\begin{equation}
f(a_{1},a_{2},a_{3}) \geq \sqrt{3a_{1}a_{2}a_{3}}.
\end{equation}
This bound is notably sharp as discussed by Beck and Zacks \cite{Beck2004}. Regarding the average behavior, Ustinov \cite{Ustinov2009} demonstrated that for three variables, the expected value of $f$ grows according to:
\begin{equation}
f(a_{1},a_{2},a_{3}) \sim \frac{8}{\pi} \sqrt{a_{1}a_{2}a_{3}}.
\end{equation}
More recent studies in additive combinatorics suggest that for a fixed $d$, the Frobenius number typically scales as 
$$
g(a_1, \dots, a_d) \sim \sqrt[d-1]{(d-1)! \prod a_i}.
$$

\subsection{The Wilf conjecture}\label{appB3}

A long-standing open question is the conjecture proposed by Wilf \cite{Wilf1978}, which relates the Frobenius number $F$, the number of omitted values (the genus $g$), and the number of generators $d$:
\begin{equation}
d \geq \frac{F+1}{F+1-g}.
\end{equation}
While the general case remains open, significant milestones have been reached. Moscariello and Sammartano \cite{Moscariello2015} validated an asymptotic version of this inequality. More recently, the conjecture has been verified for numerical semigroups with small genus ($g \leq 60$) and for those with specific patterns in their Ap\'ery sets \cite{Delgado2018}. Eliahou \cite{Eliahou2020} further refined these bounds, providing stronger evidence that the conjecture holds for semigroups with high embedding dimension, suggesting that the ratio of representable versus non-representable integers is fundamentally constrained by the density of the generating set.

\section{$\mathrm{FrobCrypt}$, a \texttt{Python} implementation}\label{appC}

The following program, named ``FrobCrypt'' enables one to test the method. FrobCrypt is a security tool that hides your messages inside a stream of numbers. To anyone else, it looks like random data. To you, it is a secret message.

\subsection{Different steps in the program}\label{appC1}

If the secret key is a set of ``magic coins'' (for example, 5 and 7). Using only these coins, some amounts are ``impossible'' to pay (like 1, 2, 3, 4, 6, 8, and 9). These are called ``gaps''. The first step is the encryption: the program takes the message and hides it inside these ``impossible'' numbers. The security is ensured by the fact that, without the secret key (the magic coins), an attacker cannot tell the difference between a normal number and a gap. The last step is the decryption: since we have the key, we can easily find the gaps and read the hidden message.

\subsection{Advantages and use of the program}\label{appC2}

The output looks like raw sensor data or random noise. The program is very fast, since it uses optimized math (Dijkstra's algorithm \cite{Dijkstra1959,Cormen2022}) to encrypt and decrypt quickly. It is also secured, since it is based on mathematical problems believed to be computationally difficult without the key. The code is in pure \texttt{Python}, and very simple to use. The procedure is as follows:

\begin{itemize}

\item Create a key and encrypt a message: from \verb|FrobCrypt| import \verb|FrobCipher|, \verb|generate_secure_key|.

\item Generate a new secret key: \verb|my_key = generate_secure_key(5)|

\verb|cipher = FrobCipher(my_key)|

\item Encrypt a word: the word 'BONJOUR' will become a list of secret numbers.

\end{itemize}

\noindent The \texttt{Python} program (with all the above steps included) reads:

\begin{verbatim}

import heapq
import random
import math

class FrobCipher:
    def __init__(self, key_set):
        """ Initializes the system with a key (set of numbers) """
        self.A = sorted(list(key_set))
        self.a1 = self.A[0]
        # Calculation of the structure of the gaps
        self.n_struct = self._build_residue_graph()
        self.F = max(self.n_struct.values()) - self.a1

    def _build_residue_graph(self):
        """ Calculate the shortest path (Dijkstra) to find the gaps """
        n = {r: float('inf') for r in range(self.a1)}
        n[0] = 0
        queue = [(0, 0)]
        while queue:
            d, r = heapq.heappop(queue)
            if d > n[r]: continue
            for x in self.A:
                nv, nr = d + x, (d + x) % self.a1
                if nv < n[nr]:
                    n[nr] = nv
                    heapq.heappush(queue, (nv, nr))
        return n

    def encrypt_byte(self, byte_val):
        """ Transforms 1 octet into 2 Frobenius numbers (4 bits each) """
        output = []
        # Splitting into nibbles (quartets)
        for nibble in [(byte_val >> 4) & 0x0F, byte_val & 0x0F]:
            found = False
            while not found:
                candidate = random.randint(1, self.F)
                # check: Is it a gap AND does it correspond to a nibble?
                if candidate < self.n_struct[candidate % self.a1] and (candidate % 16 == nibble):
                    output.append(candidate)
                    found = True
        return output

    def decrypt_byte(self, n1, n2):
        """ Recover the original octet from the two numbers """
        high = n1 % 16
        low = n2 % 16
        return (high << 4) | low

def generate_secure_key(n_elements=5):
    """ Generate a robust random key (PGCD = 1) """
    while True:
        a = random.randint(500, 1000)
        b = a + random.randint(1, 200)
        
        if math.gcd(a, b) != 1:
            continue
            
        key = [a, b]
        for _ in range(n_elements - 2):
            next_val = (random.randint(1, 3) * key[0]) + (random.randint(1, 3) * key[-1])
            if next_val not in key:
                key.append(next_val)
        
        if math.gcd(*key) == 1:
            return set(key)

# --- TEST ZONE ---
if __name__ == "__main__":
    # 1. Generation de la cle
    my_key = generate_secure_key(5)
    print(f"Key generated: {my_key}")
    
    cipher = FrobCipher(my_key)
    
    # 2. Test message
    msg = "BONJOUR"
    print(f"Original message: {msg}")

    # 3. Encryption
    encrypted_stream = []
    for char in msg:
        nums = cipher.encrypt_byte(ord(char))
        encrypted_stream.extend(nums)
    
    print(f"Encrypted stream: {encrypted_stream}")

    # 4. Decryption
    decrypted_msg = ""
    for i in range(0, len(encrypted_stream), 2):
        val = cipher.decrypt_byte(encrypted_stream[i], encrypted_stream[i+1])
        decrypted_msg += chr(val)

    print(f"Decrypted message: {decrypted_msg}")

\end{verbatim}

\end{document}